\documentclass[12pt,twoside]{article}

\usepackage{amsfonts}
\usepackage{amssymb}
\usepackage{amsmath}
\usepackage{graphicx}

\newcommand{\ignore}[1]{}

\def\@begintheorem#1#2{\par\bgroup{\sc #1\ #2. }\it\ignorespaces}
\def\@opargbegintheorem#1#2#3{\par\bgroup{\sc #1\ #2\ (#3). } \it\ignorespaces}
\def\@endtheorem{\egroup}
\newtheorem{theorem}{Theorem}[section]
\newtheorem{corollary}[theorem]{Corollary}
\newtheorem{lemma}[theorem]{Lemma}
\newtheorem{definition}[theorem]{Definition}
\newcommand{\bt}[1]{\begin{theorem}\label{#1}}
\newcommand{\bc}[1]{\begin{corollary}\label{#1}}
\newcommand{\bl}[1]{\begin{lemma}\label{#1}}
\newcommand{\ba}[1]{\begin{algorithm}\rm\label{#1}}
\newcommand{\bd}[1]{\begin{definition}\rm\label{#1}}
\newcommand{\bpr}{\noindent {\em Proof. }}
\newcommand{\et}{\end{theorem}}
\newcommand{\ec}{\end{corollary}}
\newcommand{\el}{\end{lemma}}
\newcommand{\ep}{\end{proposition}}
\newcommand{\ed}{\end{definition}}
\newcommand{\epr}{{\ \vbox{\hrule\hbox{%
\vrule height1.3ex\hskip0.8ex\vrule}\hrule}}\\\par}

\def\R{\mathbb{R}}
\def\Z{\mathbb{Z}}
\def \G {{\cal G}}
\def \l {\langle}
\def \r {\rangle}

\begin{document}

\title{\bf Multicommodity flow in Polynomial time}

\author{
Raymond Hemmecke
\and
Shmuel Onn
\thanks{Supported in part by a grant from ISF - the Israel Science Foundation}
\and
Robert Weismantel
}

\date{}
\maketitle

\begin{abstract}
The multicommodity flow problem is NP-hard already
for two commodities over bipartite graphs.
Nonetheless, using our recent theory of $n$-fold
integer programming and extensions developed herein,
we are able to establish the surprising polynomial time
solvability of the problem in two broad situations.
\end{abstract}

\section{Introduction}
\label{Introduction}

\begin{figure}[hbt]
\hskip-1.1cm
\includegraphics[scale=0.58]{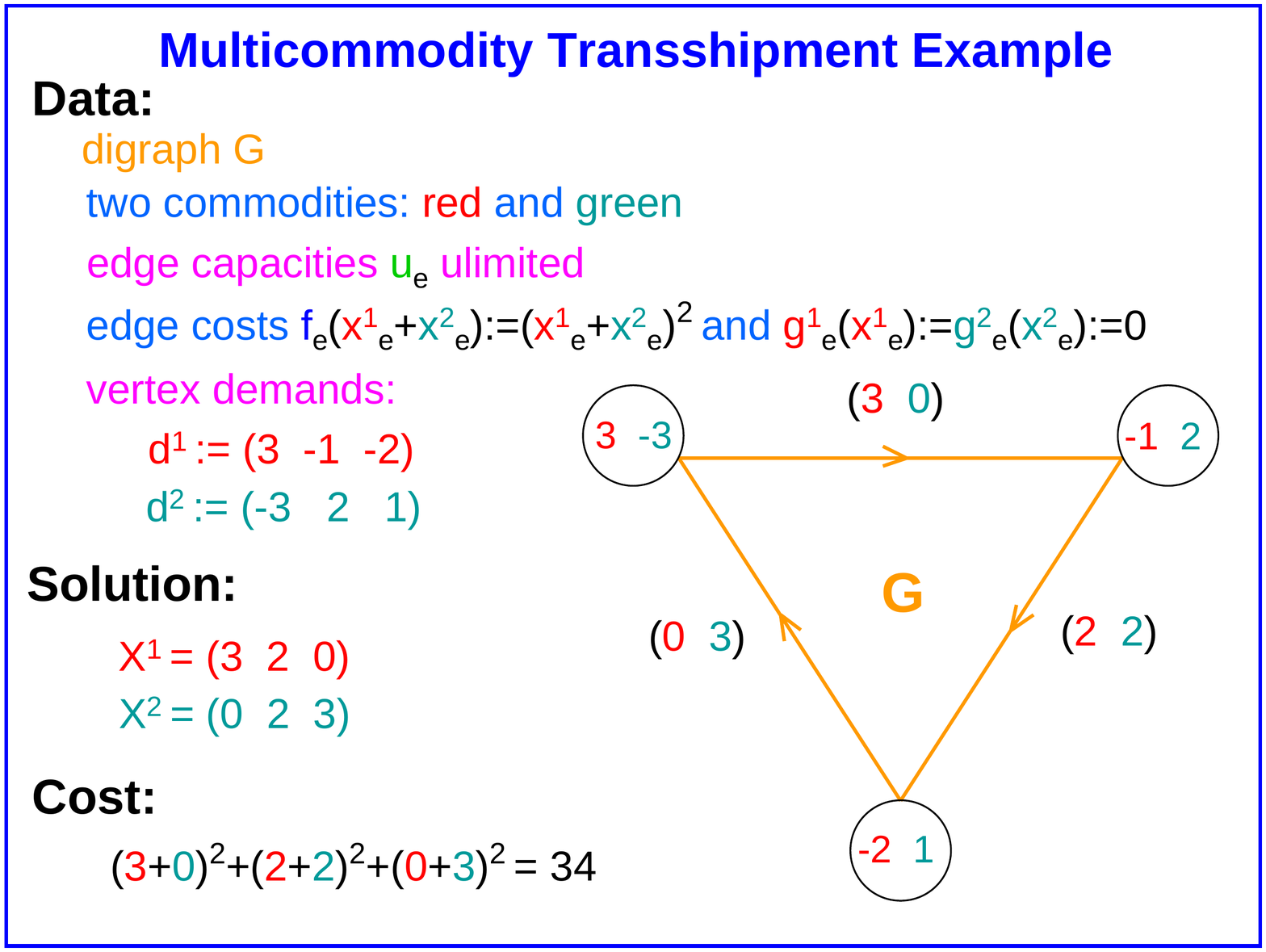}
\caption{Multicommodity Transshipment Example}
\label{multi-transshipment-figure}
\end{figure}

The multicommodity transshipment problem is a very general flow problem
which seeks minimum cost routing of several discrete commodities
over a digraph subject to vertex demand and edge capacity constraints.
The data for the problem is as follows
(see Figure \ref{multi-transshipment-figure} below for a small example).
There is a digraph $G$ with $s$ vertices and $t$ edges.
There are $l$ types of commodities.
Each commodity has a demand vector $d^k\in\Z^s$ with $d^k_v$
the demand for commodity $k$ at vertex $v$ (interpreted as supply
when positive and consumption when negative).
Each edge $e$ has a capacity $u_e$
(upper bound on the combined flow of all commodities on it).
A {\em multicommodity transshipment} is a vector $x=(x^1,\dots,x^l)$
with $x^k\in\Z_+^t$ for all $k$ and $x^k_e$
the flow of commodity $k$ on edge $e$,
satisfying the capacity constraint $\sum_{k=1}^l x^k_e\leq u_e$
for each edge $e$ and demand constraint
$\sum_{e\in\delta^+(v)}x^k_e -\sum_{e\in\delta^-(v)}x^k_e=d^k_v$
for each vertex $v$ and commodity $k$
(with $\delta^+(v),\delta^-(v)$ the sets
of edges entering and leaving vertex $v$).

The cost of transshipment $x$ is defined as follows. There are cost functions
$f_e,g^k_e:\Z\rightarrow\Z$ for each edge and each edge-commodity pair.
The transshipment cost on edge $e$
is $f_e(\sum_{k=1}^l x^k_e)+\sum_{k=1}^l g^k_e(x^k_e)$ with the first
term being the value of $f_e$ on the combined flow of all
commodities on $e$ and the second term being the sum of costs that
depend on both the edge and the commodity. The total cost is
$$\sum_{e=1}^t\left(f_e\left(\sum_{k=1}^l x^k_e\right)
+\sum_{k=1}^l g^k_e(x^k_e)\right)\ .$$

Our results apply to cost functions which can
be standard linear or convex such as
$\alpha_e|\sum_{k=1}^l x^k_e|^{\beta_e}+
\sum_{k=1}^l \gamma^k_e |x^k_e|^{\delta^k_e}$
for some nonnegative integers $\alpha_e,\beta_e,\gamma^k_e,\delta^k_e$,
which take into account the increase in cost due to channel
congestion when subject to heavy traffic or communication load
(with the linear case obtained by $\beta_e=\delta^k_e$=1).

\vskip.2cm

The problem is generally hard: even deciding if a
feasible transshipment exists (regardless of its cost)
is NP-complete already in the following two very special cases:
first, with only $l=2$ commodities over the complete bipartite
digraphs $K_{m,n}$ (oriented from one side to the other) \cite{DO1,DO2};
and second, with variable number of commodities over the
digraphs $K_{3,n}$ with $m=3$ vertices on one side (see Section 4).

Nonetheless, using the theory of $n$-fold integer programming recently
introduced in \cite{DHOW,DHORW,HOW} and extensions developed herein,
we are able to establish the surprising polynomial time solvability
of the problem, with either standard linear costs or more general costs
with {\em nonlinear convex} functions $f_e$, $g^k_e$, in two situations as follows.

\vskip.2cm
First, over any fixed digraph, we can solve the problem with a
variable number $l$ of commodities (hence termed the {\em many-commodity}
transshipment problem). This problem may seem at a first glance very
restricted: however, even for the single tiny bipartite digraph $K_{3,3}$,
we are not aware as of yet of any solution method other than the one
provided herein; and as noted, the problem is NP-hard for the digraphs $K_{3,n}$.
Our first theorem is the following (see Section 3 for the precise statement).

\bt{Transshipment}
For any fixed digraph $G$, the (convex) many-commodity transshipment problem
with variable $l$ commodities over $G$ can be solved in polynomial time.
\et
We also point out the following immediate corollary of Theorem \ref{Transshipment}.
\bc{Corollary}
For any fixed $s$, the (convex) many-commodity transshipment problem with variable
$l$ commodities on any $s$-vertex digraph is polynomial time solvable.
\ec
The complexity of the algorithm of Theorem \ref{Transshipment}
involves a term of $O\left(l^{g(G)}\right)$ where $g(G)$ is the
{\em Graver complexity} of $G$, a fascinating new digraph invariant
about which very little is known (even $g(K_{3,4})$ is as yet unknown),
see discussion in Section 4.

\vskip.2cm
Second, when the number $l$ of commodities is fixed, we can solve
the problem over any bipartite subdigraph of $K_{m,n}$
(the so-called multicommodity {\em transportation} problem)
with fixed number $m$ of suppliers and variable number $n$ of consumers.
This is very natural in operations research applications where
few facilities serve many customers. Here each commodity type $k$
may have its own volume $v_k$ per unit. Note again that if $l$
is variable then the problem is NP-hard already for $m=3$, so our
following second theorem is best possible (see Section 3 for the precise statement).

\bt{Transportation}
For fixed $l$ commodities and $m$ suppliers, the (convex)
multicommodity transportation problem with variable $n$
consumers is polynomial time solvable.
\et

We point out that the running time of our algorithms depends naturally on
the {\em binary-encoding} length $\l d^k_v,u_e\r$ of the numerical part of the
data consisting of the demands and capacities (see Section 3), so our
algorithms can handle very large numbers. To get such polynomial running
time {\em even} in the much more limited situation when {\em both} the digraph
{\em and} the number of commodities are fixed (where the number $lt$
of variables becomes fixed) and where the cost functions are linear,
one needs off-hand the algorithm of integer programming in
fixed dimension \cite{Len}. However, Theorems \ref{Transshipment} and
\ref{Transportation} involve variable dimension and \cite{Len} does not apply.

In Section 2 we review the recent theory of $n$-fold integer programming
and establish a new theorem enabling the solvability of a generalized
class of $n$-fold integer programs. In Section 3 we use the results of
Section 2 to obtain our multicommodity flow Theorems \ref{Transshipment}
and \ref{Transportation}. We conclude in Section 4 with a short discussion.

\section{$N$-fold integer programming}
\label{$N$-fold integer programming}

\subsection{Background}
\label{Background}

Linear integer programming is the following fundamental optimization problem,
$$\min\,\left\{wx\ :\ x\in\Z^n\,,\ Ax=b\,,\ l\leq x\leq u\right\}\ ,$$
where $A$ is an integer $m\times n$ matrix, $b\in\Z^m$, and $l,u\in\Z_{\infty}^n$
with $\Z_{\infty}:=\Z\uplus\{\pm\infty\}$.
It is generally NP-hard, but polynomial time solvable in two
fundamental situations: the dimension is fixed \cite{Len};
the underlying matrix is totally unimodular \cite{HK}.

Recently, in \cite{DHOW}, a new fundamental polynomial time solvable
situation was discovered. We proceed to describe this class
of so-termed {\em $n$-fold integer programs}.

An {\em $(r,s)\times t$ bimatrix}
is a matrix $A$ consisting of two blocks $A_1$, $A_2$, with $A_1$
its $r\times t$ submatrix consisting of the first $r$ rows and $A_2$
its $s\times t$ submatrix consisting of the last $s$ rows. The
{\em $n$-fold product} of $A$ is the following $(r+ns)\times nt$ matrix,
$$A^{(n)}\quad:=\quad
\left(
\begin{array}{cccc}
  A_1    & A_1    & \cdots & A_1    \\
  A_2    & 0      & \cdots & 0      \\
  0      & A_2    & \cdots & 0      \\
  \vdots & \vdots & \ddots & \vdots \\
  0      & 0      & \cdots & A_2    \\
\end{array}
\right)\quad .
$$
The following result of \cite{DHOW} asserts
that $n$-fold integer programs are efficiently solvable.
\bt{NFold1}{\bf \cite{DHOW}}
For every fixed integer $(r,s)\times t$ bimatrix $A$, there is
an algorithm that, given positive integer $n$,
$w\in\Z^{nt}$, $b\in\Z^{r+ns}$, and $l,u\in\Z_{\infty}^{nt}$,
solves in time which is polynomial in $n$ and in the binary-encoding length
$\l w,b,l,u\r$ of the rest of the data, the following so-termed
linear $n$-fold integer programming problem,
$$\min\,\left\{wx\ :\ x\in\Z^{nt}\,,\ A^{(n)} x=b\,,\ l\leq x\leq u\right\}\ .$$
\et

Some explanatory notes are in order. First, the dimension
of an $n$-fold integer program is $nt$ and is variable.
Second, $n$-fold products $A^{(n)}$ are highly non totally unimodular:
the $n$-fold product of the simple $(0,1)\times 1$ bimatrix
with $A_1$ empty and $A_2:=2$ satisfies $A^{(n)}=2 I_n$ and has exponential
determinant $2^n$. So this is indeed a class of programs which cannot
be solved by methods of fixed dimension or totally unimodular matrices.
Third, this class of programs turns out to be very natural and has
numerous applications, the most generic being to integer optimization
over multidimensional tables. In fact it is {\em universal}: the results
of \cite{DO2} imply that {\em any} integer program {\em is} an $n$-fold
program over some simple bimatrix $A$, see Section \ref{Discussion}.

The above theorem extends to $n$-fold integer programming with
nonlinear objective functions as well. The following two results,
from \cite{DHORW} and \cite{HOW} respectively, assert that the
maximization and minimization of certain convex functions over
$n$-fold integer programs can also be done in polynomial time.
The function $f$ is presented by a {\em comparison oracle} that
for any two vectors $x,y$ can check if $f(x)\leq f(y)$.

\bt{NFold2}{\bf \cite{DHORW}}
For every fixed $d$ and $(r,s)\times t$ integer bimatrix $A$,
there is an algorithm that, given $n$, bounds
$l,u\in\Z_{\infty}^{nt}$, integer $d\times {nt}$ matrix $W$,
$b\in\Z^{r+ns}$, and convex function $f:\Z^d\rightarrow\R$ presented
by a comparison oracle, solves in time polynomial in $n$ and $\l l,u,W,b \r$,
the convex $n$-fold integer maximization problem
$$\max\{f(Wx)\ :\ x\in\Z^{nt}\,,\ A^{(n)}x=b\,,\ l\leq x\leq u\}\ .$$
\et

In the next theorem, $f$ is {\em separable convex}, namely
$f(x)=\sum_i f_i(x_i)$ with each $f_i$ univariate convex.
The running time depends also on $\log {\hat f}$ with
${\hat f}$ the maximum value of $|f(x)|$ over the feasible set
(${\hat f}$ is not needed to be part of the input).

\bt{NFold3}{\bf \cite{HOW}}
For every fixed integer $(r,s)\times t$ bimatrix $A$, there is an algorithm
that, given $n$, lower and upper bounds $l,u\in\Z_{\infty}^{nt}$, $b\in\Z^{r+ns}$,
and separable convex function $f:\Z^{nt}\rightarrow\Z$
presented by a comparison oracle, solves in time which is polynomial in $n$ and
$\l l,u,b,{\hat f}\r$ the convex $n$-fold integer minimization problem
\begin{equation*}
\min\left\{f(x)\ :\ x\in\Z^{nt}\,,\ A^{(n)} x=b\,,\ l\leq x\leq u\right\}\ .
\end{equation*}
\et

\subsection{Generalization}
\label{Generalization}

We now provide a broad generalization of Theorem \ref{NFold3}
which will be useful for the multicommodity flow applications
to follow and is interesting on its own right.

We need to review some material from \cite{DHOW,HOW}.
We make use of a partial order $\sqsubseteq$ on $\R^n$ defined as follows.
For two vectors $x,y\in\R^n$ we write $x\sqsubseteq y$
if $x_iy_i\geq 0$ and $|x_i|\leq |y_i|$ for $i=1,\ldots,n$, that is, $x$ and $y$
lie in the same orthant of $\R^n$ and each component of $x$ is bounded by
the corresponding component of $y$ in absolute value. A classical lemma of Gordan
\cite{Gor} implies that every subset of $\Z^n$ has finitely-many $\sqsubseteq$-minimal
elements. The following fundamental object was introduced in \cite{Gra}.

\bd{GraverBasisDefinition}
The {\em Graver basis} of an integer matrix $A$ is defined to be the
finite set $\G(A)\subset\Z^n$ of $\sqsubseteq$-minimal elements
in $\{x\in\Z^n\,:\, Ax=0,\ x\neq 0\}$.
\ed
The Graver basis is typically exponential and cannot
be written down, let alone computed, in polynomial time.
However, we have the following lemma from \cite{DHOW}.
\bl{GraverComputation}
For every fixed integer bimatrix $A$ there is an algorithm
that, given $n$, obtains the Graver basis $\G(A^{(n)})$
of the $n$-fold product of $A$ in time polynomial in $n$.
\el
We also need the following lemma from \cite{HOW}
showing the usefulness of Graver bases.
\bl{Augmentation}
There is an algorithm that, given an integer $m\times n$ matrix $A$,
its Graver basis $\G(A)$, $l,u\in\Z_{\infty}^n$, $b\in\Z^m$,
and separable convex function $f:\Z^n\rightarrow\Z$ presented
by a comparison oracle, solves in time polynomial in
$\l A,\G(A),l,u,b,{\hat f}\r$, the program
$$\min\{f(x)\ :\ x\in\Z^n\,,\ Ax=b\,,\ l\leq x\leq u\}\ .$$
\el
Note that Lemmas \ref{GraverComputation} and \ref{Augmentation}
together imply at once Theorem \ref{NFold3} mentioned above.

\vskip.2cm
We proceed with two new lemmas needed in the proof of our generalized theorem.
\bl{ExtendedGraverComputation}
For every fixed integer $(r,s)\times t$ bimatrix $A$ and $(p,q)\times t$
bimatrix $W$, there is an algorithm that, given any positive integer $n$,
computes in time polynomial in $n$, the Graver basis $\G(B)$ of
the following $(r+ns+p+nq)\times(nt+p+nq)$ matrix,
$$B\ :=\ \left(
\begin{array}{cc}
  A^{(n)}  &  0  \\
  W^{(n)}  & I  \\
\end{array}
\right)\ .$$
\el
\bpr
Let $D$ be the $(r+p,s+q)\times (t+p+q)$ bimatrix whose blocks are defined by
$$
D_1\ :=\
\left(
\begin{array}{ccc}
A_1 & 0   & 0 \\
W_1 & I_p & 0 \\
\end{array}
\right)\ ,\quad
D_2\ :=\
\left(
\begin{array}{ccc}
A_2 & 0 & 0   \\
W_2 & 0 & I_q \\
\end{array}
\right)\quad.
$$
Apply the algorithm of Lemma \ref{GraverComputation} and
compute in polynomial time the Graver basis $\G(D^{(n)})$ of
the $n$-fold product of $D$, which is the following matrix:
$$
D^{(n)}\ =\
{\small\left(
\begin{array}{ccc|ccc|c|ccc}
A_1 & 0   & 0 & A_1 & 0   & 0 & \cdots & A_1 & 0   & 0 \\
W_1 & I_p & 0 & W_1 & I_p & 0 & \cdots & W_1 & I_p & 0 \\
\hline
A_2 & 0 & 0   & 0   & 0   & 0 & \cdots & 0   & 0   & 0 \\
W_2 & 0 & I_q & 0   & 0   & 0 & \cdots & 0   & 0   & 0 \\
\hline
0   & 0 & 0   & A_2 & 0 & 0   & \cdots & 0   & 0   & 0 \\
0   & 0 & 0   & W_2 & 0 & I_q & \cdots & 0   & 0   & 0 \\
\hline
\vdots & \vdots & \vdots & \vdots & \vdots & \vdots & \ddots & \vdots & \vdots & \vdots\\
\hline
0   & 0 & 0   & 0   & 0 & 0   & \cdots & A_2 & 0   & 0 \\
0   & 0 & 0   & 0   & 0 & 0   & \cdots & W_2 & 0   & I_q \\
\end{array}
\right)}\ .
$$
Suitable row and column permutations applied to $D^{(n)}$ give the following matrix:
$$
C\ :=\
{\small\left(
\begin{array}{cccc|cccc|cccc}
A_1 & A_1 & \cdots & A_1 & 0   & 0   & \cdots & 0    & 0   & 0   & \cdots & 0   \\
A_2 & 0   & \cdots & 0   & 0   & 0   & \cdots & 0    & 0   & 0   & \cdots & 0   \\
0   & A_2 & \cdots & 0   & 0   & 0   & \cdots & 0    & 0   & 0   & \cdots & 0   \\
\vdots & \vdots & \ddots & \vdots & \vdots & \vdots  &
\ddots & \vdots & \vdots & \vdots & \ddots & \vdots \\
0   & 0   & \cdots & A_2 & 0   & 0   & \cdots & 0    & 0   & 0   & \cdots & 0   \\
\hline
W_1 & W_1 & \cdots & W_1 & I_p & I_p & \cdots & I_p  & 0   & 0   & \cdots & 0   \\
W_2 & 0   & \cdots & 0   & 0   & 0   & \cdots & 0    & I_q & 0   & \cdots & 0   \\
0   & W_2 & \cdots & 0   & 0   & 0   & \cdots & 0    & 0   & I_q & \cdots & 0   \\
\vdots & \vdots & \ddots & \vdots & \vdots & \vdots &
\ddots & \vdots & \vdots & \vdots & \ddots & \vdots \\
0   & 0   & \cdots & W_2 & 0   & 0   & \cdots & 0    & 0   & 0   & \cdots & I_q \\
\end{array}
\right)}\ .
$$
Obtain the Graver basis $\G(C)$ in polynomial time from $\G(D^{(n)})$
by permuting the entries of each element of the latter by the
permutation of the columns of $\G(D^{(n)})$ that is used to get $C$
(the permutation of the rows does not affect the Graver basis).

Now, note that the matrix $B$ can be obtained from $C$ by dropping all
but the first $p$ columns in the second block. Consider any element
in $\G(C)$, indexed, according to the block structure, as
$(x^1,x^2,\dots,x^n,y^1,y^2,\dots,y^n,z^1,z^2,\dots,z^n)$.
Clearly, if $y^k=0$ for $k=2,\dots,n$ then the restriction
$(x^1,x^2,\dots,x^n,y^1,z^1,z^2,\dots,z^n)$ of this element is in the
Graver basis of $B$. On the other hand, if
$(x^1,x^2,\dots,x^n,y^1,z^1,z^2,\dots,z^n)$ is any element in $\G(B)$
then its extension $(x^1,x^2,\dots,x^n,y^1,0,\dots,0,z^1,z^2,\dots,z^n)$
is clearly in $\G(C)$. So the Graver basis of $B$ can be obtained
in polynomial time by
$$\G(B)\, :=\,\left\{(x^1,\dots,x^n,y^1,z^1,\dots,z^n)\, :\,
(x^1,\dots,x^n,y^1,0,\dots,0,z^1,\dots,z^n)\in\G(C)\right\}\, .
$$
This completes the proof.
\epr

In the next lemma and theorem, as before, ${\hat f}$ and ${\hat g}$
denote the maximum values of $|f(Wx)|$ and $|g(x)|$ over the feasible set
(${\hat f}$, ${\hat g}$ do not need to be part of the input).

\bl{ExtendedAugmentation}
There is an algorithm that, given an integer $m\times n$
matrix $A$, an integer $d\times n$ matrix $W$,
$l,u\in\Z_{\infty}^n$, ${\hat l},{\hat u}\in\Z_{\infty}^d$, $b\in\Z^m$,
the Graver basis $\G(B)$ of
$$B\ :=\ \left(
\begin{array}{cc}
  A  &  0  \\
  W  & I  \\
\end{array}
\right)\ ,$$
and separable convex functions $f:\Z^d\rightarrow\Z$,
$g:\Z^n\rightarrow\Z$ presented by comparison oracles, solves in
time polynomial in $\l A,W,\G(B),l,u,{\hat l},{\hat u},b,{\hat f},{\hat g}\r$,
the program
\begin{equation*}
\min\{f(Wx)+g(x)\ :\ x\in\Z^n\,,\ Ax=b
\,,\ {\hat l}\leq Wx\leq{\hat u}\,,\ l\leq x\leq u\}\ .
\end{equation*}
\el
\bpr
Define $h:\Z^{n+d}\rightarrow\Z$ by $h(x,y):=f(-y)+g(x)$ for all
$x\in\Z^n$ and $y\in\Z^d$. Clearly, $h$ is separable convex since $f,g$ are.
Now, our problem can be rewritten as
$$\min\{h(x,y):(x,y)\in\Z^{n+d},\
\left(
\begin{array}{cc}
  A  &  0  \\
  W  &  I  \\
\end{array}
\right)
\left(
\begin{array}{c}
  x  \\
  y  \\
\end{array}
\right)=
\left(
\begin{array}{c}
  b  \\
  0  \\
\end{array}
\right),\
l\leq x\leq u,-{\hat u}\leq y\leq-{\hat l}\}\ ,$$
and the statement follows at once by applying Lemma
\ref{Augmentation} to this problem.
\epr

We can now provide our new
theorem on generalized $n$-fold integer programming.
\bt{NFold4}
For every fixed integer $(r,s)\times t$ bimatrix $A$ and integer $(p,q)\times t$
bimatrix $W$, there is an algorithm that, given $n$,
$l,u\in\Z_{\infty}^{nt}$, ${\hat l},{\hat u}\in\Z_{\infty}^{p+nq}$, $b\in\Z^{r+ns}$,
and separable convex
functions $f:\Z^{p+nq}\rightarrow\Z$, $g:\Z^{nt}\rightarrow\Z$ presented
by comparison oracles, solves in time
polynomial in $n$ and $\l l,u,{\hat l},{\hat u},b,{\hat f},{\hat g}\r$,
the generalized problem
\begin{equation*}\label{NFoldStrongSeparableConvexEquation}
\min\left\{f(W^{(n)}x)+g(x)\ :\ x\in\Z^{nt}\,,\ A^{(n)}x=b\,,\
{\hat l}\leq W^{(n)}x\leq{\hat u}\,,\ l\leq x\leq u\right\}\ .
\end{equation*}
\et
\bpr
First use the algorithm of Lemma \ref{ExtendedGraverComputation}
to compute the Graver basis $\G(B)$ of
$$B\ :=\ \left(
\begin{array}{cc}
  A^{(n)}  &  0  \\
  W^{(n)}  & I  \\
\end{array}
\right)\ .$$
Now use the algorithm of Lemma \ref{ExtendedAugmentation}
to solve the problem in polynomial time.
\epr

\section{Multicommodity flows}

\subsection{Many-commodity transshipment}
\label{Many-commodity transshipment}

We begin with our theorem on nonlinear many-commodity transshipment.
As in the previous section, ${\hat f}$, ${\hat g}$ denote the maximum
absolute values of the objective functions $f$, $g$ over the feasible set.
It is usually easy to determine an upper bound on these values from the
problem data (for instance, in the special case of linear cost
functions $f$, $g$, bounds which are polynomial in the binary-encoding
length of the costs $\alpha_e$, $\gamma^k_e$, capacities $u$,
and demands $d^k_v$, readily follow from Cramer's rule).

\vskip.2cm\noindent
{\bf Theorem \ref{Transshipment}}
{\em For every fixed digraph $G$ there is an algorithm that, given
$l$ commodity types, demand $d^k_v\in\Z$ for each commodity $k$
and vertex $v$, edge capacities $u_e\in\Z_+$, and convex costs
$f_e,g^k_e:\Z\rightarrow\Z$ presented by comparison oracles,
solves in time polynomial in $l$ and $\l d^k_v,u_e,{\hat f},{\hat g}\r$,
the many-commodity transshipment problem,
\begin{eqnarray*}
&\min & \sum_e \left(f_e\left(\sum_{k=1}^l x^k_e\right)
+\sum_{k=1}^l g^k_e(x^k_e)\right) \\
&\mbox{\rm s.t.}& x^k_e\in\Z\,,\ \
\sum_{e\in\delta^+(v)}x^k_e -\sum_{e\in\delta^-(v)}x^k_e=d^k_v\,,\ \
\sum_{k=1}^l x^k_e\leq u_e\,,\ \ x^k_e\geq 0\ \ .
\end{eqnarray*}
}

\noindent
\bpr
Assume $G$ has $s$ vertices and $t$ edges and let $D$
be its $s\times t$ vertex-edge incidence matrix.
Let $f:\Z^t\rightarrow \Z$ and $g:\Z^{lt}\rightarrow \Z$ be the
separable convex functions defined by $f(y):=\sum_{e=1}^t f_e(y_e)$
with $y_e:=\sum_{k=1}^l x^k_e$
and $g(x):=\sum_{e=1}^t\sum_{k=1}^l g^k_e(x^k_e)$.
Let $x=(x^1,\dots,x^l)$ be the vector of variables
with $x^k\in\Z^t$ the flow of commodity $k$ for each $k$.
Then the problem can be rewritten in vector form as
\begin{equation*}
\min\left\{f\left(\sum_{k=1}^l x^k\right)+g\left(x\right)\ :\ x\in\Z^{lt}\,,\
Dx^k=d^k\,,\ \sum_{k=1}^l x^k\leq u\,,\ x\geq 0\right\}\ .
\end{equation*}
We can now proceed in two ways.

First way: extend the vector of variables to $x=(x^0,x^1,\dots,x^l)$
with $x^0\in\Z^t$ representing an additional slack commodity. Then the
capacity constraints become $\sum_{k=0}^l x^k=u$ and the cost function
becomes $f(u-x_0)+g(x^1,\dots,x^l)$ which is also separable convex.
Now let $A$ be the $(t,s)\times t$ bimatrix with first block
$A_1:=I_t$ the $t\times t$ identity matrix and second block $A_2:=D$.
Let $d^0:=Du-\sum_{k=1}^l d^k$ and let $b:=(u,d^0,d^1,\dots,d^l$).
Then the problem becomes the $(l+1)$-fold integer program
\begin{equation*}
\min\left\{f\left(u-x^0\right)+g\left(x^1,\dots,x^l\right)
\ :\ x\in\Z^{(l+1)t}\,,\ A^{(l)}x=b\,,\ x\geq 0\right\}\ .
\end{equation*}
By Theorem \ref{NFold3} this program
can be solved in polynomial time as claimed.

Second way: let $A$ be the $(0,s)\times t$ bimatrix with first block $A_1$ empty
and second block $A_2:=D$. Let $W$ be the $(t,0)\times t$ bimatrix with
first block $W_1:=I_t$ the $t\times t$ identity matrix and second block $W_2$ empty.
Let $b:=(d^1,\dots, d^l)$. Then the problem
is precisely the following $l$-fold integer program,
\begin{equation*}
\min\left\{f\left(W^{(l)}x\right)+g\left(x\right)
\ :\ x\in\Z^{lt}\,,\ A^{(l)}x=b\,,\ W^{(l)}x\leq u\,,\ x\geq 0\right\}\ .
\end{equation*}
By Theorem \ref{NFold4} this program
can be solved in polynomial time as claimed.
\epr

\subsection{Multicommodity transportation}
\label{Multicommodity transportation}

We proceed with our theorem on nonlinear multicommodity transportation.
The underlying digraph is $K_{m,n}$
(with edges oriented from suppliers to consumers). The problem
over any subdigraph $G$ of $K_{m,n}$ reduces to that over $K_{m,n}$
by simply forcing $0$ capacity on all edges not present in $G$.
Note that Theorem \ref{Transshipment} implies that if $m,n$ are fixed
then the problem can be solved in polynomial time for variable
number $l$ of commodities. However, we now want to allow the number $n$
of consumers to vary and fix the number $l$ of commodities instead.
This seems to be a harder problem (with no seeming analog for
non bipartite digraphs),
and the formulation below is more delicate. Therefore it is
convenient to change the labeling of the data a little bit as follows
(see Figure \ref{multi-transportation-figure} below).
\begin{figure}[hbt]
\hskip-1.1cm
\includegraphics[scale=0.58]{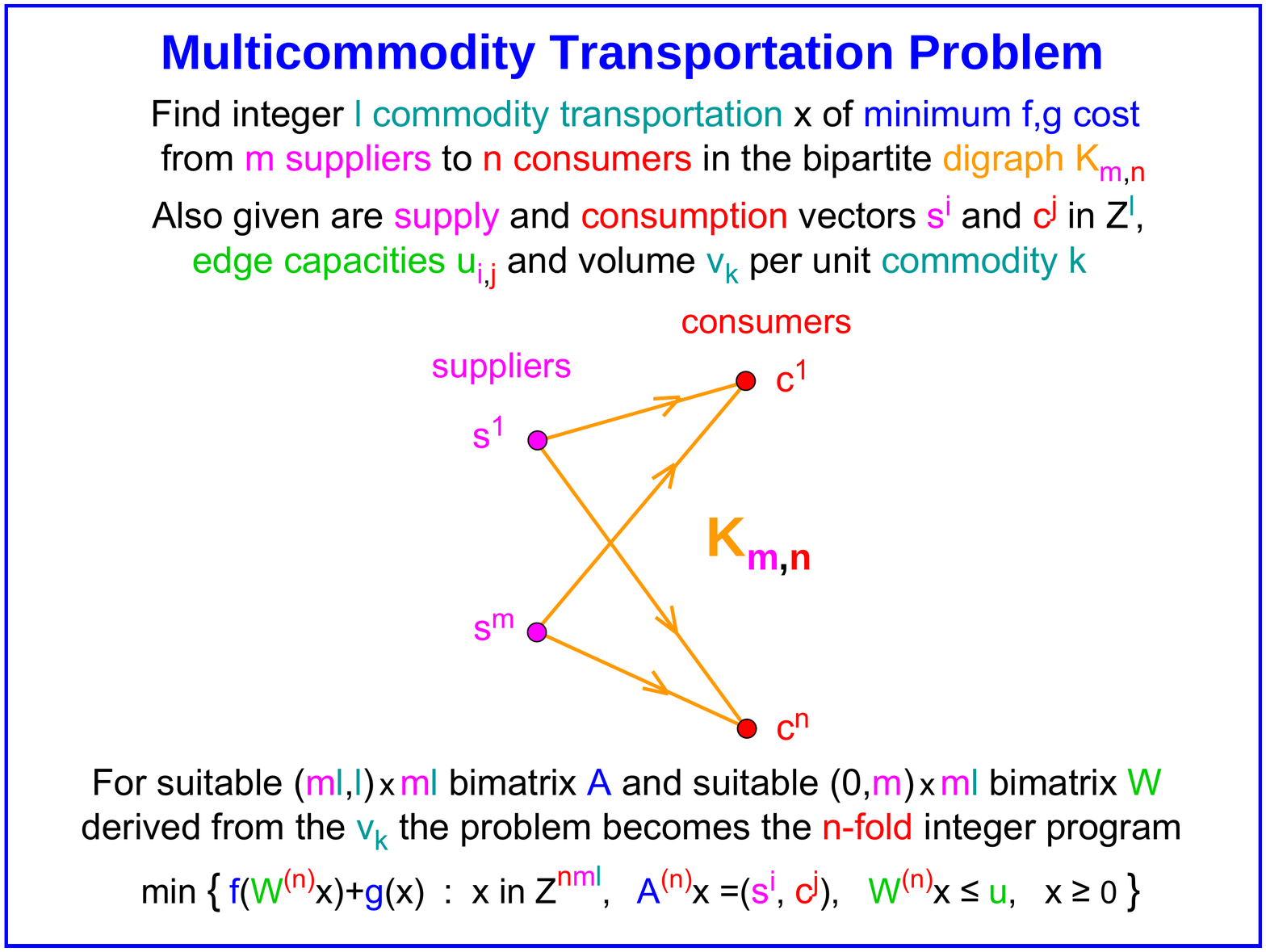}
\caption{Multicommodity Transportation Problem}
\label{multi-transportation-figure}
\end{figure}
We now denote edges by pairs $(i,j)$ where $1\leq i\leq m$ is a supplier
and $1\leq j\leq n$ is a consumer. The demand vectors are now replaced
by (nonnegative) supply and consumption vectors: each supplier $i$
has a supply vector $s^i\in\Z_+^l$ with $s^i_k$ its supply in commodity $k$,
and each consumer $j$ has a consumption vector $c^j\in\Z_+^l$ with
$c^j_k$ its consumption in commodity $k$.
In addition, each commodity $k$ may have its own volume $v_k\in\Z_+$ per unit flow.
A multicommodity transportation is now indexed as
$x=(x^1,\dots,x^n)$ with
$x^j=(x^j_{1,1},\dots,x^j_{1,l},\dots,x^j_{m,1},\dots,x^j_{m,l})$, where
$x^j_{i,k}$ is the flow of commodity $k$ from supplier $i$ to consumer $j$.
The capacity constraint on edge $(i,j)$ is
$\sum_{k=1}^l v_k x^j_{i,k}\leq u_{i,j}$ and the cost is
$f_{i,j}\left(\sum_{k=1}^l v_k x^j_{i,k}\right)+
\sum_{k=1}^l g^j_{i,k}\left(x^j_{i,k}\right)$
with $f_{i,j},g^j_{i,k}:\Z\rightarrow\Z$ convex.
As before, ${\hat f}$, ${\hat g}$ denote the maximum absolute
values of $f$, $g$ over the feasible set.

\vskip.2cm\noindent
{\bf Theorem \ref{Transportation}}
{\em For any fixed $l$ commodities, $m$ suppliers, and volumes
$v_k\in\Z_+$, there is an algorithm that, given $n$, supplies and
demands $s^i,c^j\in\Z_+^l$, capacities $u_{i,j}\in\Z_+$, and convex
costs $f_{i,j},g^j_{i,k}:\Z\rightarrow\Z$ presented by comparison oracles,
solves in time polynomial in $n$ and $\l s^i,c^j,u,{\hat f},{\hat g}\r$,
the multicommodity transportation problem,
\begin{eqnarray*}
&\min & \sum_{i,j} \left(f_{i,j}\left(\sum_k v_k x^j_{i,k}\right)
+\sum_{k=1}^l g^j_{i,k}\left(x^j_{i,k}\right)\right) \\
&\mbox{\rm s.t.}& x^j_{i,k}\in\Z\,,\ \
\sum_j x^j_{i,k}=s^i_k\,,\ \ \sum_i x^j_{i,k}=c^j_k\,,\ \
\sum_{k=1}^l v_k x^j_{i,k}\leq u_{i,j}\,,\ \ x^j_{i,k}\geq 0\ \ .
\end{eqnarray*}
}

\noindent
\bpr
Construct bimatrices $A$ and $W$ as follows.
Let $D$ be the $(l,0)\times l$ bimatrix with first block $D_1:=I_l$ and
second block $D_2$ empty. Let $V$ be the $(0,1)\times l$ bimatrix with
first block $V_1$ empty and second block $V_2:=(v_1,\dots,v_l)$.
Let $A$ be the $(ml,l)\times ml$ bimatrix with first block $A_1:=I_{ml}$
and second block $A_2:=D^{(m)}$. Let $W$ be the $(0,m)\times ml$ bimatrix
with first block $W_1$ empty and second block $W_2:=V^{(m)}$.
Let $b$ be the $(ml+nl)$-vector $b:=(s^1,\dots,s^m,c^1,\dots,c^n)$.

Let $f:\Z^{nm}\rightarrow \Z$ and $g:\Z^{nml}\rightarrow \Z$ be the
separable convex functions defined by $f(y):=\sum_{i,j}f_{i,j}(y_{i,j})$
with $y_{i,j}:=\sum_{k=1}^l v_k x^j_{i,k}$ and
$g(x):=\sum_{i,j}\sum_{k=1}^l g^j_{i,k}(x^j_{i,k})$.

Now note that $A^{(n)}x$ is an $(ml+nl)$-vector, whose first $ml$ entries
are the flows from each supplier of each commodity to all consumers,
and whose last $nl$ entries are the flows to each consumer of each
commodity from all suppliers. Therefore the supply and consumption
equations are encoded by $A^{(n)}x=b$. Next note that the
$nm$-vector $y=(y_{1,1},\dots,y_{m,1},\dots,y_{1,n},\dots,y_{m,n})$
satisfies $y=W^{(n)}x$. So the capacity constraints become $W^{(n)}x\leq u$
and the cost function becomes $f(W^{(n)}x)+g(x)$. Therefore, the
problem is precisely the following $n$-fold integer program,
\begin{equation*}
\min\left\{f\left(W^{(n)}x\right)+g\left(x\right)\ :\
x\in\Z^{nml}\,,\ A^{(n)}x=b\,,\ W^{(n)}x\leq u\,,\ x\geq 0\right\}\ .
\end{equation*}
By Theorem \ref{NFold4} this program can be solved in polynomial time as claimed.
\epr

\section{Discussion}
\label{Discussion}

We conclude with a short discussion of the universality for integer programming
of the many-commodity transportation problem and the complexity of our algorithms.

Consider the following special form of the $n$-fold product.
For an integer $s\times t$ matrix $D$, let $D^{[n]}:=A^{(n)}$
where $A$ is the $(t,s)\times t$ bimatrix $A$ with first block $A_1:=I_t$
the $t\times t$ identity matrix and second block $A_2:=D$. We consider such
$n$-fold products of the $1\times 3$ matrix ${\bf 1}_3:=[1,1,1]$.
Note that ${\bf 1}_3^{[n]}$ is precisely the $(3+n)\times 3n$
incidence matrix of the complete bipartite graph $K_{3,n}$. For instance,
\begin{equation*}\label{matrix}
{\bf 1}_3^{[3]}\ =\
\left(
\begin{array}{ccccccccc}
  1 & 0 & 0 & 1 & 0 & 0 & 1 & 0 & 0 \\
  0 & 1 & 0 & 0 & 1 & 0 & 0 & 1 & 0 \\
  0 & 0 & 1 & 0 & 0 & 1 & 0 & 0 & 1 \\
  1 & 1 & 1 & 0 & 0 & 0 & 0 & 0 & 0 \\
  0 & 0 & 0 & 1 & 1 & 1 & 0 & 0 & 0 \\
  0 & 0 & 0 & 0 & 0 & 0 & 1 & 1 & 1 \\
\end{array}
\right)\ .
\end{equation*}
The following surprising theorem was proved in \cite{DO2} building
on results of \cite{DO1}. (For further details and consequences
for privacy in statistical databases see \cite{DO2,DO3,Onn}.)

\vskip.2cm\noindent
{\bf The Universality Theorem \cite{DO2}}
{\em Every rational polytope $\{x\in\R_+^d\,:\,Bx=b\}$ stands in polynomial
time computable integer preserving bijection with some polytope
\begin{equation}\label{universality_equation}
\left\{x\in\R_+^{3nl}\ :\ {\bf 1}_3^{[n][l]}x=a\right\}\ .
\end{equation}}

\vskip.2cm\noindent
In particular, every integer program can be lifted in polynomial time to
a program over a matrix ${\bf 1}_3^{[n][l]}$ which is completely determined
by two parameters $n$ and $l$ only.

Now note (see proof of Theorem \ref{Transshipment})
that the integer points in (\ref{universality_equation})
are {\em precisely} the feasible points of some $(l-1)$-commodity transshipment
problem over $K_{3,n}$. So every integer program can be lifted in
polynomial time to some $l$-commodity program over some $K_{3,n}$.
Thus, the many-commodity transportation problem, already
over the digraphs $K_{3,n}$ with fixed number $3$ of suppliers, is {\em universal}
for integer programming. So, in particular, the $l$-commodity transportation
problem over $K_{3,n}$ is NP-hard when both $n,l$ are variable, but polynomial
time solvable for arbitrary fixed number $n$ of consumers and variable
number $l$ of commodities by Theorem \ref{Transshipment}.

\vskip.2cm
Our algorithms involve two major tasks: the construction of the Graver
basis of a suitable $n$-fold product in Lemmas \ref{GraverComputation}
and \ref{ExtendedGraverComputation}, and the iterative use of
this Graver basis to solve the underlying (convex) integer program in Lemmas
\ref{Augmentation} and \ref{ExtendedAugmentation}. The polynomial time
solvability of these tasks is established in \cite{DHOW,HOW}. Here we only briefly
discuss the complexity of the first task in the special case of a digraph, which
is relevant for the complexity of the many-commodity transshipment application.

Let $D$ be the $s\times t$ incidence matrix of a digraph $G$.
Consider $l$-fold products $D^{[l]}$ of the special form defined above.
The {\em type} of an element $x=(x^1,\dots,x^l)$ in the Graver basis
$\G(D^{[l]})$ is the number of nonzero blocks $x^k\in\Z^t$ of $x$.
It turns out that for any digraph $G$ there is a finite nonnegative
integer $g(G)$ which is the largest type of any element of any
$\G(D^{[l]})$ independent of $l$. We call this new digraph
invariant $g(G)$ the {\em Graver complexity} of $G$. The complexity
of computing $\G(D^{[l]})$ is $O(l^{g(G)})$ (see \cite{DHOW})
and hence the importance of $g(G)$. Unfortunately, our present understanding of
the Graver complexity of a digraph is very limited and much more study
is required. Very little is known even for the complete bipartite digraphs
$K_{3,n}$ (oriented from one side to the other): while $g(K_{3,3})=9$,
already $g(K_{3,4})$ is unknown. See \cite{BO} for more details and a
lower bound on $g(K_{3,n})$ which is exponential in $n$.

\section*{Acknowledgements}

The work of Shmuel Onn on this article was mostly done
while he was visiting and delivering the Nachdiplom Lectures
at ETH Z\"urich. He would like to thank Komei Fukuda
and Hans-Jakob L\"uthi for related stimulating discussions
during this period.

\noindent {\small Raymond Hemmecke}\newline
\emph{Otto-von-Guericke Universit\"at Magdeburg,
D-39106 Magdeburg, Germany}\newline
\emph{email: hemmecke{\small @}imo.math.uni-magdeburg.de}\newline
\ \ \emph{http://www.math.uni-magdeburg.de/{\small $\sim$hemmecke}}

\vskip0.5cm
\noindent {\small Shmuel Onn}\newline
\emph{Technion - Israel Institute of Technology, 32000 Haifa, Israel}\newline
and\newline
\emph{ETH Z\"urich, 8092 Z\"urich, Switzerland}\newline
\emph{email: onn{\small @}ie.technion.ac.il}\newline
\emph{http://ie.technion.ac.il/{\small $\sim$onn}}

\vskip0.5cm
\noindent {\small Robert Weismantel}\newline
\emph{Otto-von-Guericke Universit\"at Magdeburg,
D-39106 Magdeburg, Germany}\newline
\emph{email: weismantel{\small @}imo.math.uni-magdeburg.de}\newline
\emph{http://www.math.uni-magdeburg.de/{\small $\sim$weismant}}

\end{document}